\documentstyle[12pt,amsfonts,amssymb]{article}

\newcommand{\proof}{{\it Proof.$\:\:\:\:$}}

\newcommand{\taaa}{{\frak t}}
\newcommand{\haaa}{{\frak h}}

\newcommand{\R}{{\Bbb R}}

\newcommand{\C}{{\Bbb C}}

\newcommand{\gaaa}{{\frak g}}
\newcommand{\maaa}{{\frak m}}
\newcommand{\aaaa}{{\frak a}}

\newcommand{\id}{{ \mathrm{ id}}}

\newcommand{\nat}{{\Bbb  N}}

\newcommand{\aca}{{\aaaa_\C^\ast}}

\def\hB{\hspace*{\fill}$\Box$\newline\noindent}

\newtheorem{prop}{Proposition}
\newtheorem{lem}[prop]{Lemma}

\newtheorem{theorem}[prop]{Theorem}
\newtheorem{kor}[prop]{Corollary}

\def\imath{i}

\normalbaselines 
\begin{document} 
\begin{center}
\vspace*{1.0cm}

{\LARGE{\bf Hodge theory on hyperbolic manifolds of infinite volume}} 

\vskip 1.5cm

{\large {\bf Martin Olbrich }} 

\vskip 0.5 cm 

Mathematisches Institut \\ 
Universit\"at G\"ottingen \\ 
Bunsenstr. 3-5 \\ 
D-37073 G\"ottingen \\ Germany

\end{center}

\vspace{1 cm}

\begin{abstract}
Let $X=H^n$ be the real hyperbolic space of dimension $n$, $G=SO(1,n)_0$ the
identity component of its group of isometries, and let 
$\Gamma\subset G$ be a discrete torsion-free subgroup
acting convex-cocompactly on $X$. We are going to study the (de Rham) cohomology groups $H^p(Y,\C)\cong H^p(\Gamma,\C)$ of the hyperbolic manifold
$Y=\Gamma\backslash X$. In this paper we describe certain spaces
of generalized $\Gamma$-invariant currents on the sphere at infinity of $X$ 
with support
on the limit set of $\Gamma$. These spaces are finite-dimensional. The main
result of the present paper identifies the cohomology with a quotient of such
spaces (such a result was conjectured by Patterson \cite{patterson93}). We explain in which sense this result generalizes the classical Hodge Theorem for the case of compact quotients. We obtain analogous results for the groups $H^p(\Gamma,F)$, where $F$ is an irreducible
finite-dimensional representation of $G$. 
\end{abstract}

\vspace{1 cm} 

\section{Cocompact and convex-cocompact subgroups}

Let $X=H^n$ be the real hyperbolic space of dimension $n$, and let
$\partial X$ be its geodesic boundary. The union $X\cup \partial X$
can be given the structure of an analytic manifold with boundary (use, e.g.,
the unit ball model of $H^n$) with an analytic action of the group $G=SO(1,n)_0$.

We consider a discrete torsion-free subgroup $\Gamma\subset G$. Its
limit set $\Lambda\subset \partial X$ is defined to be the set of accumulation
points in $X\cup \partial X$ of the $\Gamma$-orbit of an arbitrary point $x\in X$. The complement
$\Omega:=\partial X\setminus\Lambda$ is called the domain of discontinuity
of $\Gamma$. Indeed, $\Gamma$ acts properly discontinuously on $\Omega$ as
well as on $X\cup \Omega$. Hence, $\overline Y:=\Gamma\backslash (X\cup \Omega)$
is manifold with boundary $B:= \Gamma\backslash \Omega$. Its interior
$Y:=\Gamma\backslash X$ carries the structure of a complete Riemannian manifold
with constant sectional curvature $-1$, i.e., a hyperbolic manifold.
One calls $\Gamma$ convex-cocompact if $\overline Y$ is compact. In particular,
a cocompact subgroup is convex-cocompact. In this case $\Lambda=\partial X$,
$B=\emptyset$, and $\overline Y=Y$ is a closed hyperbolic manifold.

Our first goal is to reformulate the classical Hodge theorem for closed
Riemannian manifolds in the special case of closed hyperbolic manifolds, i.e.,
cocompact subgroups $\Gamma\subset G$, in such a way that it still would make good
sense for convex-cocompact non-cocompact $\Gamma$. The main tool will
be the Poisson transform which sends differential forms on $\partial X$ to harmonic forms on the interior. Let $\Omega^p(X)$, $\Omega^p_{-\omega}(\partial X)$
denote the spaces of smooth $p$-forms on $X$ and $p$-forms with hyperfunction coefficients on $\partial X$, respectively. 
We consider the differential
$d$, the codifferential $\delta=d^*$ defined by the hyperbolic metric, and the form Laplacian $\Delta=\delta d+d\delta$. Gaillard proved the following

\begin{theorem}[\cite{Gai1}]\label{key}
There are Poisson transforms which provide $G$-equivariant isomorphisms
$$\Omega^p_{-\omega}(\partial X)\cong \{\omega\in\Omega^p(X)\:|\:\Delta\omega=0,\delta\omega=0\}\ ,\quad p\ne\frac{n-1}{2}\ ,$$
and
$$\Omega^\pm_{-\omega}(\partial X)\cong \{\omega\in\Omega^\frac{n-1}{2}(X)\:|\:d\omega=0,\delta\omega=0\}\ ,$$
where $\Omega^\pm_{-\omega}(\partial X)$ are the eigenspaces of the conformally invariant $*$-operator 
$$*:\Omega^\frac{n-1}{2}(\partial X)\rightarrow
\Omega^\frac{n-1}{2}(\partial X)\ .$$
\end{theorem}
The theorem remains true if one replaces the left hand sides by the smaller spaces $\Omega^*_{-\infty}(\partial X)$ of forms with distribution coefficients,
i.e., currents, and adds on the right hand sides the condition that the harmonic
forms should have moderate growth (see \cite{Gai1}). Nowadays one version of the theorem can be obtained from the other by a "change of globalization
functor" (see e.g. \cite{casselman89}). 

By a slight abuse of notation we denote by $\Omega^*_{-\infty}(\Lambda) \subset\Omega^*_{-\infty}(\partial X)$ the space of all currents supported on the limit set of $\Gamma$. $\Gamma$ acts on these spaces, and we denote by ${}^\Gamma\Omega^*_{-\infty}(\Lambda)$ and ${}^\Gamma\Omega^*_{-\infty}(\partial X)$ the corresponding spaces
of $\Gamma$-invariants.  
Here is the promised reformulation of the Hodge Theorem.

\begin{kor}\label{obv}  
If $Y=\Gamma\backslash X$ is a closed hyperbolic manifold, then there are the following isomorphisms
\begin{eqnarray*}
 H^p(Y,\C)&\cong & {}^\Gamma\Omega^{n-p}_{-\infty}(\Lambda)\ ,\quad p\not\in\{0,\frac{n+1}{2}\}\ ,\\
H^\frac{n+1}{2}(Y,\C)&\cong& {}^\Gamma\Omega^\pm_{-\infty}(\Lambda)\ .
\end{eqnarray*}
\end{kor}
\proof Recall that in the cocompact case $\Lambda=\partial X$. Thus by Theorem \ref{key}
and the compactness of $Y$ (in particular, moderate growth is automatic)
$${}^\Gamma\Omega^*_{-\infty}(\Lambda)\cong \{\omega\in\Omega^*(Y)\:|\:d\omega=0,\delta\omega=0\}{\cong}\{\omega\in\Omega^{n-*}(Y)\:|\:d\omega=0,\delta\omega=0\}\ .$$
The second isomorphism is given by the $*$-operator. The right hand side is isomorphic
to $H^{n-*}(Y,\C)$ by the Hodge Theorem.
\hB
The unmotivated incorporation of the $*$-operator in the isomorphism is necessary in
order to produce a statement which has some chance to generalize to convex-cocompact $\Gamma$. Indeed, it turns out that Corollary \ref{obv} remains true in the general case for
$p\ge \frac{n+1}{2}$. For $p\le\frac{n}{2}$ sometimes a modification is necessary, see Theorem \ref{main} below.

\section{The extension map}\label{ext}

From now on let $\Gamma\subset G$ be convex-cocompact. In this section we assume in addition that
$\Omega\ne\emptyset$, i.e., $\Gamma$ is not cocompact.

The aim of this section is to review the theory of the extension operator $ext$ which
extends $\Gamma$-invariant hyperfunction sections defined on $\Omega$ of homogeneous vector bundles over $\partial X$ across the limit set in order to get globally defined ones.
This theory was developed by U. Bunke and the author in \cite{bunkeolbrich963} and \cite{bunkeolbrich982}. It will enable us to get hold on the spaces  ${}^\Gamma\Omega^{p}_{-\infty}(\Lambda)$ and its relatives defined below.

As a homogeneous space we have $\partial X=G/P$, where $P=MAN$ is the group of
orientation preserving conformal transformations of $\R^{n-1}$ leaving the point $\infty$ fixed. Here $M=SO(n-1)$, $A\cong \R_+$ is the group of dilatations, and $N$ is the group
of translations. Let $\aaaa$ be the Lie algebra of $A$. We will identify its complexified dual $\aca$ with $\C$. Let $\hat M$ be the set of equivalence classes of irreducible representations of $M$ in complex vector spaces. Any pair $(\sigma,\lambda)\in\hat M\times \aca$ gives rise to
an irreducible representation $\sigma_\lambda$ of $P$ in the representation space of $\sigma$ by 
$$\sigma_\lambda(man):=\sigma(m) a^{\rho-\lambda}\ ,\quad \rho=\frac{n-1}{2}\ .$$
Let $V(\sigma_\lambda)\rightarrow\partial X$ be the associated homogeneous vector bundle.
By $C^{-\omega}(\partial X, V(\sigma_\lambda))$ we denote the space of its hyperfunction sections. It carries a representation of $G$ known under the name principal series representation. Let $\sigma^p$ be the $p$-th exterior power of the standard representation of $SO(n-1)$. Then we have
$$\Omega^p_{-\omega}(\partial X)=C^{-\omega}(\partial X, V(\sigma^p_{\rho-p}))\ .$$

We want to understand the space of $\Gamma$-invariant sections ${}^\Gamma C^{-\omega}(\partial X, V(\sigma_\lambda))$, in particular its subspace ${}^\Gamma C^{-\omega}(\Lambda, V(\sigma_\lambda))$ of sections
supported on the limit set. Since $\Gamma\backslash \Omega$ is a compact manifold the space of invariant
sections along $\Omega$ is easy to understand:
$$ {}^\Gamma C^{-\omega}(\Omega, V(\sigma_\lambda))\cong C^{-\omega}(B, V_B(\sigma_\lambda))\ .$$
It is the space of hyperfunction sections of the corresponding bundle $V_B(\sigma(\lambda))$ on $B$. Restriction of a hyperfunction on $\partial X$ to the
open set $\Omega$ provides a map
$$ res:{}^\Gamma C^{-\omega}(\partial X, V(\sigma_\lambda))\rightarrow 
C^{-\omega}(B, V_B(\sigma_\lambda))\ .$$
It turns out that for generic $\lambda\in\aca$ the map $res$ is an isomorphism. 
Moreover, $res$ has an inverse which depends meromorphically on $\lambda$:

\begin{theorem}[\cite{bunkeolbrich982}, \cite{bunkeolbrich963}]\label{BO}
There is a meromorphic family of maps with finite-dimensional singularities
$$ext_\lambda: C^{-\omega}(B, V_B(\sigma_\lambda))\rightarrow {}^\Gamma C^{-\omega}(\partial X, V(\sigma_\lambda))$$
such that 
\begin{equation}\label{mds}
res\circ ext_\lambda=\id\ .
\end{equation}
\end{theorem}

A singularity at $\lambda$ is called finite-dimensional if the coefficients in the
principal part of the Laurent expansion at $\lambda$ are finite-dimensional operators.
The theorem is proved in \cite{bunkeolbrich982} in the distribution
setting, only. But the arguments of \cite{bunkeolbrich963}, 2.5. and 2.6. show how
one can transfer it to hyperfunctions.

In order to deal with points $\lambda\in\aca$, where $ext$ has a pole, we have
to introduce slightly more general bundles. We consider the representation
$1^k: MAN\rightarrow \C^k$ given by the Jordan block
$$ 1^k(man)=\exp\left(\begin{array}{ccccc}
0&\log(a)&0&\dots& 0\\
0&0&\log(a)&\dots&0\\
&&\ddots&&\\
0&0&\dots&0&\log(a)\\
0&0&\dots&0&0 \end{array}\right)\ .$$
Then we set $V^k(\sigma_\lambda):=V(\sigma_\lambda\otimes 1^k)$. Obviously, $V^1(\sigma_\lambda)=V(\sigma_\lambda)$, and for any $l\in\nat$ there is a short exact sequence
$$ 0\rightarrow V^k(\sigma_\lambda)\stackrel{i_l}{\longrightarrow}
V^{k+l}(\sigma_\lambda)\stackrel{p_k}{\longrightarrow} V^l(\sigma_\lambda)\rightarrow 0\ .$$
Again we have a restriction map between the corresponding spaces of sections
$$ res:{}^\Gamma C^{-\omega}(\partial X, V^k(\sigma_\lambda))\rightarrow 
C^{-\omega}(B, V^k_B(\sigma_\lambda))\ .$$ 
The space of sections of the extended bundle $C^{-\omega}(\partial X, V^k(\sigma_\lambda))$ can be $G$-equivariantly identified with the quotient of the space of germs at $\lambda$ of meromorphic families $\mu\mapsto f_\mu\in C^{-\omega}(\partial X, V(\sigma_\mu))$ which have a pole of order at most $k$ by the space of germs of
holomorphic families. By choosing a holomorphic identification of the bundles
$V_B(\sigma_\mu)$ for varying $\mu$ with $V_B(\sigma_\lambda)$ we can think of
$V^k_B(\sigma_\lambda)$ as a subspace of the germs at $\lambda$ of meromorphic families $\mu\mapsto f_\mu\in C^{-\omega}(B, V_B(\sigma_\mu))$.  These interpretations show that for any $k\in\nat$ the family of operators $ext_\mu$ induces a map
$$ ext^k: C^{-\omega}(B, V^k_B(\sigma_\lambda))\rightarrow {}^\Gamma C^{-\omega}(\partial X, V^{k+k_-}(\sigma_\lambda))\ ,$$
where
$k_-=k_-(\sigma_\lambda)$ denotes the order of the pole of $ext$ at $\lambda$. By (\ref{mds}) it satisfies $res\circ ext^k= i_{k_-}$.  Since $res$ commutes
with $p_k$ we obtain 
\begin{equation}\label{tor}
res\circ p_k \circ ext^k=0\ .
\end{equation}

We are interested in the spaces ${}^\Gamma C^{-\infty}(\Lambda, V^{k}(\sigma_\lambda))$
of $\Gamma$-invariant distribution sections of $V^{k}(\sigma_\lambda)$ supported on $\Lambda$ which generalize the spaces of currents ${}^\Gamma\Omega^{p}_{-\infty}(\Lambda)$ discussed before.

Showing that Poisson transforms of elements of ${}^\Gamma C^{-\omega}(\Lambda, V^{k}(\sigma_\lambda))$ are always of moderate growth on $X$ one obtains 

\begin{lem}
${}^\Gamma C^{-\infty}(\Lambda, V^{k}(\sigma_\lambda))={}^\Gamma C^{-\omega}(\Lambda, V^{k}(\sigma_\lambda))$.
\end{lem}

Combining the lemma with Equation (\ref{tor}) we see that  $p_k \circ ext^k$
maps $C^{-\omega}(B, V^k_B(\sigma_\lambda))$ to ${}^\Gamma C^{-\infty}(\Lambda, V^{k_-}(\sigma_\lambda))$. We denote the image of $p_{k_-} \circ ext^{k_-}$ by $E^+_\Lambda(\sigma_\lambda)$. In fact, it consists
of all invariant distributions on the limit set which can be constructed by means of the extension map.   

Along the lines of the proofs of Proposition 6.11 and Corollary
6.12 in \cite{bunkeolbrich982} one shows in addition

\begin{prop}\label{gg}
The sequence of inclusions 
$$ {}^\Gamma C^{-\infty}(\Lambda, V^{1}(\sigma_\lambda))\subset\dots\subset {}^\Gamma C^{-\infty}(\Lambda, V^{k}(\sigma_\lambda))\subset{}^\Gamma C^{-\infty}(\Lambda, V^{k+1}(\sigma_\lambda))\subset \dots                             $$
stablizes at some $k=:k_+(\sigma_\lambda)$, and $\dim{}^\Gamma C^{-\infty}(\Lambda, V^{k_+}(\sigma_\lambda))<\infty$.
\end{prop}

\section{The main theorem}

We now want to relate the spaces ${}^\Gamma C^{-\infty}(\Lambda, V^k(\sigma^{n-p}_{p-1-\rho}))$ to cohomology. This can be done by
the following proposition which is again a corollary of Gaillard's results 
\cite{Gai1}.

\begin{prop}\label{bor}
Suitably defined Poisson transforms $P^k_p$, $P^k_\pm$ provide $G$-equivariant isomorphisms
$$P^k_p: C^{-\omega}(\partial X, V^k(\sigma^{n-p}_{p-1-\rho}))\stackrel{\cong }\longrightarrow\{\omega\in\Omega^p(X)\:|\:\Delta^k\omega=0, d\omega=0\}\ ,\quad p\ne\frac{n+1}{2}\ ,$$
and
$$P^k_\pm: C^{-\omega}(\partial X, V^k(\sigma^{\pm}_0))\stackrel{\cong }\longrightarrow \{\omega\in\Omega^\frac{n+1}{2}(X)\:|\: (d*)^k\omega=0, d\omega=0\}\ .$$
\end{prop} 

Note that the Poisson transforms $P^1_p$ differ from the ones in Theorem \ref{key}
by the application of the $*$-operator. We are now able to state the main theorem of the paper which generalizes 
Corollary \ref{obv} to arbitrary convex-cocompact groups. 

\begin{theorem}\label{main}
For $p\in \{1,2,\dots, n\}$ set $l(p):=k_+(\sigma^{n-p}_{p-1-\rho})$. Then $l(p)\le k_-(\sigma^{n-p}_{p-1-\rho})+1$, and for $k\ge l(p)$ the Poisson transforms $P^k_p$, $P^k_\pm$ induce isomorphisms
\begin{eqnarray*} 
H^p(Y,\C)&\cong& {}^\Gamma C^{-\infty}(\Lambda, V^k(\sigma^{n-p}_{p-1-\rho}))/E^+_\Lambda(\sigma^{n-p}_{p-1-\rho})\ ,\quad 
p\not=\frac{n+1}{2}\ ,\\
H^\frac{n+1}{2}(Y,\C)&\cong& {}^\Gamma C^{-\infty}(\Lambda, V^k(\sigma^{\pm}_0))/E^+_\Lambda(\sigma^\pm_0)\ .
\end{eqnarray*}
For $p\ge\frac{n+1}{2}$ we have $l(p)\le 1$ and the following isomorphisms
\begin{eqnarray}
H^p(Y,\C)&\cong & \left\{\begin{array}{clc}
{}^\Gamma\Omega^{n-p}_{-\infty}(\Lambda)\ ,&p>\frac{n+1}{2}\\    
{}^\Gamma\Omega^\pm_{-\infty}(\Lambda)\ ,&p=\frac{n+1}{2} \end{array}\right.
\label{ei}\\
&\cong& \left\{\omega\in {}^\Gamma\Omega^{n-p}_{-\infty}(\Lambda)\:|\: d\omega=0\right\}\ .
\label{eiei}
\end{eqnarray}
\end{theorem}
Let us add a couple of remarks concerning the theorem.
\begin{itemize}
\item If $\Gamma$ is cocompact then it is natural to set $E^+_\Lambda(\sigma^{n-p}_{p-1-\rho}):=\{0\}$ and $k_-(\sigma^{n-p}_{p-1-\rho})=0$. With these conventions Theorem \ref{main} applied to cocompact $\Gamma$
is just a reformulation of Corollary \ref{obv}.
\item For $p>\frac{n+1}{2}$ the theorem is a fairly direct consequence
of the results of Mazzeo/ Philipps \cite{mazzeophillips90} concerning the
$L^2$-cohomology of $Y$. This was already noted by Patterson \cite{patterson93} and was also observed
by Lott \cite{lott98}. But for $p<\frac{n+1}{2}$
$L^2$-methods are not sufficient.
\item That the cohomology of $Y$ (as well as $H^p(\Gamma,F)$ considered in the next section) should be representable by currents on
the limit set was conjectured by Patterson \cite{patterson93}. In fact,
he suggested that the isomorphism (\ref{eiei}) should be true for all 
$p$. But there are examples where the map $\left\{\omega\in {}^\Gamma\Omega^{n-p}_{-\infty}(\Lambda)\:|\: d\omega=0\right\}\rightarrow
H^p(Y,\C)$ is not injective. We do not know whether this map is surjective
in general.
\end{itemize}

We will now describe the main ideas entering
the proof of Theorem \ref{main}. Full details and several related results and applications will
appear in a forthcoming paper. 

Using the surjectivity of the Laplacian
on $\Omega^p(Y)$ one shows that any cohomology class can be represented
by a closed form which is annihilated by some power of the Laplace operator
(in fact even by a coclosed harmonic form).
Then the idea was that such a form is exact if and only if it arises as a value or a derivative at $\lambda=0$
of a family of closed eigenforms $\omega_\lambda$, $\Delta \omega_\lambda=\lambda \omega_\lambda$. Using the theory of Poisson transforms
one translates the problem to the boundary. Families of eigenforms
correspond to families of $\Gamma$-invariant hyperfunction sections of the corresponding bundles
over $\partial X$. Such families can always be described using $ext_\lambda$. Now one uses the extension
operator in order to produce a representative of each cohomology class which is supported on the limit
set. 

From the technical point of view it is simpler to replace this family
criterion by the following infinitesimal version: a closed form is exact
if and only if it belongs to the image of $\Delta^k$ restricted to
closed forms for some or, equivalently, for all $k$. Applying Proposition \ref{bor} and observing that the
action of the Laplacian on closed forms annihilated by $\Delta^k$ corresponds
to the action of $p_1$ on ${}^\Gamma C^{-\omega}(\partial X, V^k(\sigma^{n-p}_{p-1-\rho}))$ we are lead to the setting of Section \ref{ext}.
Using Proposition \ref{gg} one shows that an element is in the image of
$p_1$ if and only if it is in the image of some $ext^l$ or $p_m\circ ext^l$. Eventually we end up
with Theorem \ref{main}.

\section{Finite-dimensional $G$-representations as coefficients}

We can interpret the cohomology of $Y$ as the group cohomology
of $\Gamma$ with coefficients in its trivial representation:
$H^p(Y,\C)\cong H^p(\Gamma,\C)$. If $F$ is a finite-dimensional representation
of $G$, then we can look at it as a representation of $\Gamma$ and form
$H^p(\Gamma, F)$. Since any finite-dimensional $G$-representation is
semisimple we can assume that $F$ is irreducible. It turns out that the theory
explained in the previous two sections is well-behaved with respect to the
translation functor, i.e., tensoring with $F$
followed by the projection according to the generalized infinitesimal character of $F$.
This allows one to translate Theorem \ref{main} to $H^*(\Gamma,F)$ for
any irreducible $F$. In order to formulate the corresponding theorem
we need some preparation.

Denote by $\gaaa$ and $\maaa$ the Lie algebras of $G$ and $M$, respectively.
Let $\taaa\subset \maaa$ be the Lie algebra of a maximal torus of $M$. Then $\haaa:=\aaaa\oplus\taaa$ is a Cartan subalgebra of $\gaaa$. 
We fix a positive
Weyl chamber $\taaa_+^*\subset i\taaa^*$ and denote by $\rho_\maaa\in\taaa^*_+$ the half sum of the corresponding positive roots of $\maaa$.
As a positive Weyl chamber $\haaa^*_+ \subset \haaa^*_\R:=\aaaa^*\oplus i\taaa^*$ we choose the
one which contains $\aaaa^*_+$ and is contained in $\aaaa^*_+\oplus \taaa_+^*$.
We consider the following subset of
the Weyl group $W(\gaaa_\C,\haaa_\C)$ of $\gaaa_\C$:
$$ W^1:=\{w\in W(\gaaa_\C,\haaa_\C)\:|\: \alpha\in\haaa_+^* \Rightarrow w\alpha_{|\taaa}\in\taaa_+^*\}\ . $$
If $n$ is even then $W^1=\{w_0,w_1,\dots,w_{n-1}\}$, where $w_i$ is the unique element of length $i$ in $W^1$. For odd $n$ we have $W^1=\{w_0,\dots, w_{\frac{n-3}{2}},w_+,w_-,w_\frac{n+1}{2},\dots,w_{n-1}\}$, where $w_i$ and $w_\pm$ have length $i$ and $\frac{n-1}{2}$, respectively.

For a finite-dimensional irreducible representation $F$ with highest weight 
$\nu\in\haaa_+^*$
we define $\sigma^p_F\in \hat M$ by its highest weight
$$ \mu_p:=w_{n-1-p}(\nu+\rho+\rho_\maaa)_{|\taaa}-\rho_\maaa\in\taaa^*_+\ .$$
For odd $n$ we can replace $w_p$ by $w_\mp$ and obtain the highest weights of 
representations $\sigma^\pm_F\in \hat M$. Set $\sigma^\frac{n-1}{2}_{F}:=
\sigma^+_F\oplus\sigma^-_F$. Furthermore we define
$$\lambda_p:=-w_{n-1-p}(\nu+\rho+\rho_\maaa)_{|\aaaa}\in\aaaa^*\ .$$
The analogously defined elements $\lambda_\pm$ are equal to $0$.

Note that for $F=\C$, i.e., $\nu=0$, we have $\sigma_{\C,\lambda_p}^p=\sigma^p_{\rho-p}$. For general $F$ the parameters
$\sigma^p_{F,\lambda_p}$ explicitly describe the effect of the translation functor applied to the de Rham complex of $\partial X$. Indeed, the spaces $C^{-\omega}(\partial X, V(\sigma^{p}_{F,\lambda_{p}}))$ fit into an exact de Rham like complex of
$G$-equivariant differential operators 
\begin{eqnarray*}
0\rightarrow F\rightarrow C^{-\omega}(\partial X, V(\sigma^{0}_{F,\lambda_{0}}))&\stackrel{d_F^0}\longrightarrow & C^{-\omega}(\partial X, V(\sigma^{1}_{F,\lambda_{1}}))\stackrel{d_F^1}\longrightarrow \\
\dots&\stackrel{d_F^{n-2}}\longrightarrow& C^{-\omega}(\partial X, V(\sigma^{n-1}_{F,\lambda_{n-1}}))\rightarrow F\rightarrow 0
\end{eqnarray*}
which appears in the literature under various names like BGG-resolution or
\v Zelobenko complex. Generalizations of it also played a central role
in J. Slov\'ak's talk at the Workshop.
The following
theorem is the natural generalization of Theorem \ref{main}.

\begin{theorem}\label{weyl}
For a finite-dimensional irreducible representation $F$ of $G$ and $p\in \{1,2,\dots, n\}$ set $l(F,p):=k_+(\sigma^{n-p}_{F,\lambda_{n-p}})$. Then $l(p)\le k_-(\sigma^{n-p}_{F,\lambda_{n-p}})+1$, and for $k\ge l(F,p)$ there are isomorphisms
\begin{eqnarray*} 
H^p(\Gamma,F)&\cong& {}^\Gamma C^{-\infty}(\Lambda, V^k(\sigma^{n-p}_{F,\lambda_{n-p}}))/E^+_\Lambda(\sigma^{n-p}_{F,\lambda_{n-p}})\ ,\quad 
p\not=\frac{n+1}{2}\ ,\\
H^\frac{n+1}{2}(\Gamma,F)&\cong& {}^\Gamma C^{-\infty}(\Lambda, V^k(\sigma^{\pm}_{F,0}))/E^+_\Lambda(\sigma^\pm_{F,0})\ .
\end{eqnarray*}
For $p\ge\frac{n+1}{2}$ we have $l(p)\le 1$ and the following isomorphisms
\begin{eqnarray*}
H^p(\Gamma,F)&\cong&  \left\{\begin{array}{clc}
{}^\Gamma C^{-\infty}(\Lambda, V(\sigma^{n-p}_{F,\lambda_{n-p}}))\ ,&p>\frac{n+1}{2}\\    
{}^\Gamma C^{-\infty}(\Lambda, V(\sigma^{\pm}_{F,0}))\ ,&p=\frac{n+1}{2} \end{array}\right.\\
&\cong& \left\{\omega\in{}^\Gamma C^{-\infty}(\Lambda, V(\sigma^{n-p}_{F,\lambda_{n-p}}))\:|\: d^{n-p}_F\omega=0\right\}\ .
\end{eqnarray*}
All these isomorphisms are induced by Poisson transforms into
closed forms with values in the flat bundle over $Y$ corresponding
to $F$.
\end{theorem}

\section*{Acknowledgements} 
Discussions with the following mathematicians have strongly influenced 
my changing view on the subject over the years: A. Juhl, S. J. Patterson,
U. Bunke and P.-Y. Gaillard. Let me thank all of them. In addition, I want
to thank Joachim Hilgert who invited me to present the above results at the   
Workshop "Lie Theory and Its Applications in Physics - Lie III".

\end{document}